\input amstex
\documentstyle{amsppt} 
% \magnification=1200

%\pagewidth{40 pc}

%%%%%%%%%%%%%%%%%%%%% User defined commands %%%%%%%%%%%%%%%%%%%%%%%%%%%

%\input commands.tex

% Caligraphic Symbols:
\def\cal{\Cal}
\define\Y{\Cal Y}

\define\M{\Cal M}
\define\V{\Cal V}

\define\F{\Cal F}
\define\W{\Cal W}
\define\CalS{\Cal S}

% Bold Face Symbols:
\define\R{{\Bbb R}}
\define\Q{{\Bbb Q}}
\define\Z{{\Bbb Z}}
\define\C{{\Bbb C}}

% Greek Letters:
\define\a{\alpha}
\undefine\b
\define\b{\beta}
\define\g{\gamma}
\undefine\l
\define\l{\lambda}

% Composite Greek Symbols:
\define\lam{\Lambda^{-1,-1}}
\define\nn{\eta_{-}}
\define\np{\eta_{+}}
\define\nz{\eta_0}
\define\pip{\pi_{+}}
\define\pin{\pi_{-}}
\define\piz{\pi_0}
\define\pil{\pi_{\Lambda}}

% Lie Groups/Algebras:
\define\Gc{G_{\Bbb C}}

% \define\gg{}				% How to make a german g in tex?

% Operators:

% Misc:

\define\half{\frac{1}{2}}
\define\hph{\hphantom}
\define\pd{\partial}
\define\rel{{}^r}
\define\prel{{}^{(r)}}

%%%%%%%%%%%%%%%%%%%%%%%%% End Commands File %%%%%%%%%%%%%%%%%%%%%%%%%%%%%%%%

\define\itemb{\item"$\bullet$"}

\TagsOnRight
\TagsAsMath
\loadmsbm

\topmatter

\abstract We give a condition for a variation of mixed Hodge structure 
on a curve to be admissible. It involves  the asymptotic
behavior of a grading of the weight filtration, supplementing exactly 
the description of the graded variation and its monodromy given by Schmid's Orbit Theorems. In many salient cases the condition is equivalent to admissibility. 
\endabstract

\author Aroldo Kaplan, \\  Gregory Pearlstein
\endauthor

\leftheadtext{Aroldo Kaplan,\quad Gregory Pearlstein}

\title Singularities of variations of \\ mixed Hodge structure \endtitle

\thanks We wish to thank P.Deligne, whose ideas on limits of mixed 
Hodge structures were generously shared and play an important role here.  The significance of the grading $\cal Y$  was discussed in 
correspondence of his with E. Cattani and one of us (Kaplan) dating back to 
1992 \cite{6}.  Related results were presented in \cite{10}, \cite{12} and \cite{14}. 
\endthanks

\endtopmatter

\document

% Be sure to uncomment this before producing final draft.
\NoBlackBoxes

%%%%%%%%%%%%%%%%%%%%%%%%%%%%%% section-1 %%%%%%%%%%%%%%%%%%%%%%%%%%%%%%%%

$$
	\text{\bf 1.\quad Introduction}
$$

\par By the work of Griffiths and Schmid, the Gauss--Manin connection of a 
variation of pure, polarized Hodge structure $\V\to\Delta^*$  with unipotent 
monodromy has a regular singular point at the origin.  Accordingly \cite{5}, 
the Hodge filtration $\F$ of such a variation extends to a system of 
holomorphic subbundles of the canonical extension $\tilde\V\rightarrow 
\Delta$ of $\V$  
\cite{8}.  Moreover, as a consequence of the $SL_2$ Orbit Theorem 
\cite{15}, the limiting Hodge filtration of such a variation pairs
with the (shifted) monodromy weight filtration of $\V$ to define a limiting
mixed Hodge structure on the central fiber of $\tilde\V$.
%\vskip 10pt

\par In contrast, the situation for variations of graded-polarized mixed
Hodge structure is markedly different \cite{16}: irregular singularities and 
monodromies occur in the simplest of non-geometric examples.  For variations
with unipotent monodromy, such irregularities can be tamed by
imposing the following {\it admissibility} conditions:
\smallskip
(1) the limiting Hodge filtration  $F_{\infty}$ of $\V$ exists;

(2) the relative weight filtration $\rel W = \rel W(N,W)$ exists;
\smallskip
\noindent which in turn imply:
\smallskip
(3) the pair $(F_{\infty},\rel W)$ is a mixed Hodge structure, relative to which the 
      monodromy logarithm $N$ is morphism of type $(-1,-1)$.
\smallskip
\par
Admissible variations are closed under degenerations. Most notably, variations of 
pure structure are automatically admissible: this is a major consequence of Schmid's Orbit Theorems.  
Yet, this very fact limits their use to situations where  the Orbit Theorems can be ignored {\it ab initio}, that is,  when the limiting structure can be constructed from other data, as in \cite{16}, where it comes from geometry. It is clear that a mixed version of Griffith's program and applications like those of \cite{2},\cite{4} and \cite{11}, require a condition {\it supplementary} to Schmid's Theorems to determine admissibility. For elementary reasons, these theorems do hold for the graded variation, hence the condition must be on the degeneration of the extension data.

\par
 In this article we propose a condition -that certain smooth grading $\Bbb C$-grading of the weight filtration extends continuously over the puncture- and prove its equivalence with admissibility in a number of cases. This equivalence is a mixed analog of the SL(2)-Orbit Theorem, since the way the logarithmic monodromy breaks up under the limiting grading  characterizes ``the monodromies that can occur". 

\vskip 10pt

\par 
To state the results in detail, recall that a mixed Hodge structure $(F,W)$ on a complex vector
space $V=V_{\R}\otimes\C$ defines a unique bigrading 
$$
        V = \bigoplus_{p,q}\, I^{p,q}                           
$$
of $V$ with the following  properties:
$$F^p = \oplus_{a\geq p,b}\, I^{a,b}\qquad\qquad W_k = \oplus_{a+b\leq k}\, I^{a,b}$$
$$
	\bar I^{p,q} \equiv I^{q,p} \ \  \mod \bigoplus_{r<q,s<p}\, I^{r,s} \tag{1.1}$$

\par In particular, a mixed Hodge structure
determines a grading $Y_{(F,W)}$ of the underlying
weight filtration $W$ by the rule
$$
        Y_{(F,W)}(v) = kv \iff v\in\bigoplus_{p+q=k}\, I^{p,q}       \tag{1.2}
$$
These identification of a $\Bbb Z$-grading of $V$ with the semisimple endomorphism with degrees as eigenvalues, will be used throughout. Under it, the action of $GL(V)$ on gradings becomes the adjoint action on $ \text{End}(V)$.

\par
Our first result specifies the condition and asserts the equivalence in the case of variations that are unipotent in the sense of [9].

\proclaim{Theorem I} Let $\V\to\Delta^*$ be a variation of graded-polarized
mixed Hodge structure such that the
graded variation is constant.  Then, $\V$ is admissible if
and only if the $C^{\infty}$-grading $\Y$ of $\W$ defined pointwise
by $(1.2)$ extends continuously to  $\tilde\V\to\Delta$.
\endproclaim

\par
The implications in 
$$\V \ \text{admissible} \iff \Y 
\ \text{extends}$$
can be established in many other cases. Here, we also establish the full equivalence when the limiting objects are real in the sense specified below. 
\par
In spite of the existing evidence, it is not obvious that this reality assumption can be dropped altogether. One cannot just replace a variation of mixed structure by a split one, as one does in the pure case: in the terms of [12], the monopole $g(z)$ become singular at infinity. 
In any case, the real situation arises often enough in practice (e.g., mirror symmetry) and is a necessary step towards the general case -because it is so for pure structures.

\par
Recall that a mixed Hodge structure $(F,\rel W)$ is 
split over $\R$ if and only if 
$$
\overline{I^{p,q}_{(F,\rel W)}} = I^{q,p}_{(F,\rel W)}
$$
Accordingly, we shall say that an admissible nilpotent orbit 
$$(e^{zN}.F,W)$$
of mixed Hodge structure is {\it split} if the mixed Hodge structure 
$$(F,\rel W(N,W))$$ splits over $\Bbb R$. Note that if a nilpotent orbit $(e^{zN}.F,W)$ is admissible, that is precisely its limiting mixed Hodge structure.
Our second result is:

\proclaim{Theorem II} Let 
$$z\to e^{zN}.F$$
be a nilpotent orbit of graded-polarized mixed Hodge
structure such that the graded orbits are split. Then  
$e^{zN}.F$  is admissible and split if and only if
the limit
$$
 	Y_{\infty} = \lim_{Im(z)\to\infty}\, Ad(e^{-zN}) Y_{(e^{zN}.F,W)}
$$
exists and is defined over $\Bbb R$.
\endproclaim

\par
For a general variation, the limiting mixed Hodge structure is only defined once a coordinate on the base -more specifically, an element of $T_0(\Delta)^*$, has been chosen. We can then refer to split admissible variations, as those whose limiting mixed Hodge structure is split. With this understood, one has

\proclaim{Theorem III} Let $\V\to\Delta^*$ be a variation of 
graded-polarized mixed Hodge structure with unipotent monodromy, and 
$\Y$ denote the $C^{\infty}$ grading the weight filtration $\W$ obtained by 
applying  $(1.2)$ to each fiber of $\V$. Suppose that the graded variations are split. Then, the variation $\V$ 
is admissible and split  if and 
only if
 $\Y$ extends continuously to a grading of $\W$ in  $\tilde \V$, with limiting value $\Y(0)=\lim_{s\to 0}\, \Y(s)$ defined 
over $\R$.
\endproclaim

In any of the above situations, let
$$N = N_0 + N_{-2} + \cdots$$
be the decomposition of the monodromy logarithm relative to the grading $ Y_{\infty}$ (i.e.,  $\Y(0)$). Let 
$$
	\rho:sl_2(\C)\to End(V)		 		
$$
be the representation determined by the Orbit Theorems for pure structures applied to $Gr^W$, together with the identification $Gr^W\approx V$ provided by $Y_{\infty}$.
Let
$$
   n_-   = \pmatrix 
  	0 & 0 \\
 	1 & 0 \\
        \endpmatrix,\qquad
   h   = \pmatrix 
  	1 & 0  \\
  	0 & -1 \\ 
        \endpmatrix,\qquad
   n_+   = \pmatrix
  	0 & 1 \\
  	0 & 0 \\
         \endpmatrix				  $$
In the situation of Theorems II and III
where $Y_{\infty}$ is real, $\rho$ is defined over $\R$:
$$
	\rho:sl_2(\R)\to gl(V_{\R}).				 $$
In Theorem I instead, $\rho$ is trivial in $Gr^W$ and $N=N_{-2}$, but $Y_{\infty}$ need not be real. In any case, the following holds 

\proclaim{Theorem IV} 
{\roster
\itemb $N_0=\rho(n_-)$, 
\itemb $N_{-1}=0$
\endroster}
and for all $k>1$, $N_{-k}$ is a highest weight vector for $\rho$ of weight $k-2$, i.e., 
{\roster
\itemb $\rho(h)N_{-k}=(k-2)N_{-k}$
\itemb $\rho(n_+)N_{-k}=0.$
\endroster}
\endproclaim

\par 
Next, we sketch the proofs of these theorems. Theorem I is proved in \S 3.  One first writes the period map $F:U\to\M$ as
$$
	F(z) = e^{zN}e^{\Gamma(s)}.F_{\infty}
$$
where $F_{\infty}$ is the limiting Hodge filtration of $\V$, and 
$\Gamma(s)$ is a holomorphic function of $s=e^{2\pi iz}$ which takes 
values in the nilpotent Lie algebra 
$$
	Lie_{-1} = \{\, \a\in End(V) \mid \a:W_k\to W_{k-1}\,\forall k\,\}
$$
and vanishes at $s=0$. To finish the proof we apply Theorem $(2.9)$,
which provides a description of how the decomposition $V = \oplus_{p,q}\, I^{p,q}$ varies under the action
of $\exp(Lie_{-1})$ on $\M$.

\

\par Theorems II and III  are proved in \S 4 and \S 5, respectively. One of the implications in the first depends heavily upon the results
of Deligne \cite{6} discussed in the appendix at the end of our paper. They imply that if $W$ and $N$ arise from a nilpotent orbit of graded-polarized mixed Hodge structure $e^{zN}.F$ which is admissible and split, then
$$
	Y_{(e^{zN}.F,W)} = Ad(e^{zN})\tilde Y		
$$
where
$\tilde Y = Y(N,\rel Y)$, with $\rel Y = Y_{(F,\rel W)}$, is a grading of $W$ defined over $\R$.
In particular, whenever $e^{zN}.F$ is a split admissible nilpotent orbit, the associated grading
$$
	e^{-zN}.Y_{(e^{zN}.F,W)}
$$ 
in  Theorem II has constant value 
$Y_{\infty} = \tilde Y$.  Likewise, if $(F,\rel W)$ is split over $\R$ 
then so are each of the induced mixed Hodge structures in the graded, and the corresponding nilpotent orbits are split.

\par To prove the converse, we recall that split nilpotent orbits of mixed Hodge structure are just $SL(2,R)$-orbits, assemble the corresponding representations in the graded into one on $V$ as already described, define 
$$
	\prel Y := Y_{\infty} + \rho(h), \qquad \prel W_k = \bigoplus_{j\leq k}\, E_j(\prel Y)  $$
and prove that $\prel W_k$ is actually the relative weight filtration of $N$ and $W$. This is the most technically complicated part of the proof, involving the use of the $sl_2$ symmetry to deduce 
properties of certain Laurent series.  Theorem IV is a corollary of these proofs.

\smallskip

Some of the material from [6] and [10] is being published here for the first time.

% END SECTION 1
$$
	\text{\bf 2.\quad Preliminaries}
$$

\par In this section we reformulate our theorems in terms of
period maps, review the definition of admissibility, discuss the geometry of the associated classifying
spaces  and describe an
analog of the Nilpotent Orbit Theorem for 
variations of graded-polarized mixed Hodge structure.

\definition{Definition} A variation of graded-polarized mixed 
Hodge structure consists of a $\Q$-local system $\V_{\Q}\to S$ endowed
with the following additional data:
\roster
\item A rational, increasing weight filtration 
$$
        0\subseteq\cdots\subseteq\W_{k-1}\subseteq\W_k
         \subseteq\cdots\subseteq\V_{\C}
$$
of $\V_{\C}=\V_{\Q}\otimes\C$.
\item A decreasing Hodge filtration 
$$
        0\subseteq\cdots\subseteq\F^p\subseteq\F^{p-1}
         \subseteq\cdots\subseteq\V_{\C}\otimes {\Cal O}_S
$$
which is holomorphic and horizontal with respect to the Gauss-Manin connection
$\nabla$ of $\V_{\C}\otimes{\Cal O}_S$ 
(i.e. $\nabla\F^p \subseteq \Omega^1_S\otimes\F^{p-1}$) and pairs with the
weight filtration $\W$ to define a mixed Hodge structure on each fiber of 
$\V$.
\item A collection of rational, non-degenerate bilinear forms
$$
        \CalS_k:Gr^{\W}_k(\V_{\Q})\otimes Gr^{\W}_k(\V_{\Q})\to\Q,\qquad 
        Gr^{\W}_k := \W_k/\W_{k-1}
$$
of alternating parity $(-1)^k$ which polarize the corresponding variations
of pure Hodge structure 
$$
        \F Gr^{\W}_k := \frac{\F^p\cap\W_k + \W_{k-1}}{\W_{k-1}}
$$
\endroster
\enddefinition

Actually,
we shall make no reference to either the 
graded-polarizations or the rational structure 
itself, so
the results hold for arbitrary real variations 
of graded-polarizable mixed Hodge structure.
\vskip 10pt

\par Now, as discussed in \cite{10}, \cite{14} and elsewhere, the data of 
such a variation $\V\to S$ may be reformulated in terms of the monodromy
representation
$$
        \rho:\pi_1(S,s_o)\to Aut(\V_{s_o}),\qquad\text{Image}(\rho) = \Gamma
$$
of $\V$ on a given fiber $\V_{s_o}$ and the corresponding period map
$$
        \varphi:S\to\M/\Gamma                                
$$
obtained by parallel translating the data of $\V$ to $\V_{s_o}$.
Alternatively, upon passage to the universal cover $\tilde S$ of $S$,  such a  variation is equivalent to the monodromy 
representation $\rho$ defined above together with a $\pi_1$-equivariant map
$$
        F:\tilde S\to\M
$$
from $\tilde S$ into a suitable classifying space $\M$ of graded-polarized 
mixed Hodge structure. $F$ is both holomorphic and horizontal: i.e.
$$
        \frac{\pd}{\pd \bar z_j} F^p(z) \subseteq F^p(z),\qquad
        \frac{\pd}{\pd z_j} F^p(z) \subseteq F^{p-1}(z)
$$

In particular, if $\V\to\Delta^*$ is a variation of graded-polarized
mixed Hodge structure with unipotent monodromy action $T = e^{-N}$, its period map $\varphi$ may be 
viewed as a holomorphic, horizontal map
$$
        F: U\to\M
$$
from the upper half-plane $U$ into $\M$ which satisfies the quasi-periodicity
condition 
$$F(z+1) = e^N.F(z).$$
Our choice of sign for the monodromy logarithm
reflects the fact that we are pulling back the data of $\V$ to a fixed
reference fiber.

\par More specifically, parallel 
translation under  $\nabla$ endows our fixed 
reference fiber $V = \V_{s_0}$ with a choice of rational structure $V_{\Q}$  
as well as a constant, rational weight filtration $W$ and a collection of 
constant, rational, non-degenerate bilinear forms
$$
        \CalS_k:Gr^W_k(V_{\Q})\otimes Gr^W(V_{\Q})\to \Q
$$
of alternating parity $(-1)^k$.  Modulo the action of the 
monodromy group $\Gamma$, we also obtain a  Hodge filtration $F_s$,
which pairs with $W$ to define a graded-polarized mixed Hodge structure 
with constant graded Hodge numbers 
$$
        h^{p,q} =  \text{dim}\, H^{p,q}_s
$$
\flushpar Consider the corresponding classifying space
$$
        \M = \M(W,\CalS,h^{p,q}),
$$
consisting of all filtrations $F$ of $V$ such that $(F,W)$ is a mixed Hodge structure, which is graded-polarized by 
$\CalS$ and such that 
$$
        \text{dim}_{\C}\, F^p Gr^W_k = \sum_{r\geq p}\, h^{r,k-r}.
$$

\par As shown in \cite{10} and \cite{14}, $\M$ is a homogeneous complex
manifold which fits into an ascending sequence of homogeneous spaces
$$
        M_{\R}\subseteq\M\subseteq\check\M\subseteq\check\F(W)
              \subseteq\check\F
$$
defined as follows:
\roster
\itemb $\check\F$ is the flag variety consisting of all 
decreasing filtrations $F$ such that 
$\text{dim}_{\C}\, F^p = \sum_{r\geq p}\, h^{r,s}$. 

\itemb $\check\F(W)$ is the submanifold of $\check\F$ consisting 
of those filtrations $F\in\check\F$ such that
$\text{dim}_{\C}\, F^p Gr^W_k = \sum_{r\geq p}\, h^{r,k-r}$.

\itemb $\check\M$ is the submanifold of $\check\F(W)$ consisting
of all filtrations $F\in\check\F(W)$ which satisfy Riemann's first bilinear
relation with respect to the graded-polarizations $\CalS$.

\itemb $\M_{\R}$ is the $C^{\infty}$-submanifold of $\M$ consisting
of the filtration $F\in\M$ for which the associated mixed Hodge structure 
$(F,W)$ is split over $\R$.
\endroster
\flushpar The corresponding sequence of Lie groups is
$$
        G_{\R}\subseteq G \subseteq \Gc \subseteq GL(V)^W \subseteq GL(V)
$$
where
\roster
\itemb $GL(V)^W = \{\, g\in GL(V) \mid g:W_k\to W_k\quad\forall k\,\}$. 
\itemb $\Gc = \{\, g\in GL(V)^W \mid Gr(g) \in Aut_{\C}(\CalS)\,\}$.
\itemb $G = \{\, g\in\Gc\mid Gr(g) \in Aut_{\R}(\CalS)\,\}$.
\itemb $G_{\R} = \{\, g\in G \mid g\in GL(V_{\R})\,\}$.
\endroster
i.e. $GL(V)$ acts transitively on $\check\Cal F$, $GL(V)^W$ acts transitively
on $\check\F(W)$, $\Gc$ acts transitively on $\check\M$, $G$ acts transitively
on $\M$ and $G_{\R}$ acts transitively on $\M_{\R}$.

\

\par Next we describe the
canonical extension of $\V\to\Delta^*$ in more detail \cite{8}. Given a flat vector bundle $E\to\Delta^{*n}$ with unipotent
monodromy, there exists a unique extension $E^c\to\Delta^n$ relative to
which the flat connection of $E$ has at worst simple poles with nilpotent
residues along the divisor $D=\Delta^n/\Delta^{*n}$.  Alternatively, given 
a choice of local coordinates $(s_1,\dots,s_n)$ on $\Delta^n$ relative to 
which the divisor $D$ assumes the form $s_1\cdots s_n = 0$, the canonical 
extension $E^c$ described above may be identified with the local free sheaf 
generated by the sections
$$
        \tilde \sigma 
        = \exp(\frac{1}{2\pi i}\sum_{j=1}^n\, \log s_j\, N_j)\sigma
$$
where
$\sigma$ is a flat, multivalued section of $E$ and 
$$N_j = -\log(T_j)$$ is the monodromy logarithm 
associated to pulling back along the $j$'th loop 
$$\g_j(t) = (s_1,\dots,e^{2\pi i t}s_j,\dots, s_n).$$

\

\flushpar{\it Remark}. We shall denote
the adjoint action of $GL(V)$ upon $End(V)$ by  $g.\a$ and the action of $\Gc$ [resp. $G$] upon  $\check\M$
[resp. $\M$] by $g.F$.
\vskip 10pt

\proclaim{Lemma 2.1} Let $\V\to\Delta^*$ be a variation of graded-polarized 
mixed Hodge structure with unipotent monodromy action $T = e^{-N}$, and 
$F:U\to\M$ be a lift of the associated period map 
$\varphi:\Delta^*\to\M/\Gamma$ to the upper half-plane.   Then, the 
functions
$$
	\psi(z) = e^{-zN}.F(z),\qquad 
	Y(z) = e^{-zN}.Y_{(F(z),W)}
$$
satisfy the periodicity conditions 
$$
	\psi(z+1) = \psi(z),\qquad Y(z+1) = Y(z)
$$
and hence descend to well defined functions $\psi(s)$ and $Y(s)$ on the 
punctured disk $\Delta^*$ via the covering map $z\to s=e^{2\pi i z}$.  Moreover,
\roster
\itemb The Hodge filtration $\F$ of $\V$ extends to the canonical extension $\tilde \V\rightarrow \Delta$ if and only if the limiting Hodge filtration
$$
	F_{\infty} = \lim_{s\to 0}\, \psi(s)
$$
exists as an element of $\check\M$.

\itemb The grading $\Y$ of the underlying weight filtration $\W$ defined by
the $I^{p,q}$'s of $\V$ extends continuously to the canonical extension
of $\V$ as a grading of $\W$ if and only if 
$$
	Y_{\infty} = \lim_{s\to 0}\, Y(s)
$$
exists, and grades $\W$.
\endroster
\endproclaim
\demo{Proof} One simply identifies  $\tilde \V$ with 
the locally free sheaf generated by the sections $\tilde\sigma$ described
above.  
\enddemo

For a coordinate-free description, both the limiting Hodge 
filtration $F_{\infty}$ and the limiting grading $Y_{\infty}$ defined above 
should actually be viewed as objects defined on the cotangent space of 
$\Delta$ at zero.  More precisely, if $\tilde s$ and $s$ are holomorphic
coordiantes on $\Delta$ which vanish at $0$, then corresponding limiting
objects $F_{\infty}$ and $Y_{\infty}$ defined above will in general agree 
if and only if $(d\tilde s)_0 = (ds)_0$.
\vskip 10pt

\par Recall that a variation
of graded-polarized mixed Hodge structure with unipotent monodromy action 
$T = e^{-N}$ is {\it admissible} \cite{16} if and only if the following hold:
\roster
\item The Hodge filtration $\F$ of $\V\to\Delta^*$ extends 
holomorphically to the canonical extension $\tilde\V$.
\item The relative weight filtration $\rel W = \rel W(N,W)$ exists.
\endroster

The relative weight filtration is defined by the Lemma below. In the appendix to
\cite{16}, Deligne proved that whenever a variation of graded-polarized mixed Hodge 
structure is admissible, then the limiting Hodge 
filtration $F_{\infty}$ of $\V$ pairs with the corresponding relative weight 
filtration $\rel W$ to define a mixed Hodge structure for which $N$ is a 
$(-1,-1)$-morphism.

\vskip 10pt

\par To define the filtration $\rel W = \rel W(N,W)$, recall first that given a nilpotent endomorphism $N$ of a finite
dimensional complex vector space $V$, the corresponding
(monodromy) weight filtration $W(N)$ is the unique increasing filtration of $V$ such that
$\ell$:
\roster
\itemb $N:W_{\ell}(N) \to W_{\ell-2}(N)$.
\itemb The induced map 
$
	N^{\ell}:Gr^{W(N)}_{\ell}\to Gr^{W(N)}_{-\ell}
$
is an isomorphism.
\endroster

Given an increasing filtration $W$ of a finite 
dimensional vector space $V$ and an index $\ell\in\Z$, the corresponding
shifted object is the filtration 
$$
	W[\ell]_j = W_{j+\ell}					
$$

\proclaim{Lemma 2.2} (cf. \cite{16}) Let $W$ be an increasing filtration of a 
finite dimensional vector space $V$ and $N$ be a nilpotent endomorphism 
of $V$ which preserves $W$.  Then, there exists at most one increasing 
filtration 
$\rel W = \rel W(N,W)$ of $V$ such that for each index $k$:
\roster
\item"(i)"  $N:\rel W_k\to\rel W_{k-2}$.
\item"(ii)" $\rel W$ induces on $Gr^W_k$ the shifted monodromy weight 
filtration
$$
	W(N:Gr^W_k\to Gr^W_k)[-k]
$$
of $N$ on $Gr^W$.
\endroster
\endproclaim

\par To close this section we discuss some of the geometry of the classifying
space $\M$.  

% Was (3.2) in old draft.
\proclaim{Lemma 2.3} Let $F(z)$ be the period map of a variation of 
graded-polarized mixed Hodge structure with unipotent monodromy for which 
the limiting Hodge filtration $F_{\infty} = \lim_{Im(z)\to\infty} e^{-zN}.F(z)$
exists.  Then, given a choice of a vector space decomposition
$$
	Lie(\Gc) = Lie(\Gc^{F_{\infty}})\oplus q
$$
there exists a unique holomorphic function
$\Gamma:\Delta\to q$, $\Gamma(0) = 0$,  
$s=e^{2\pi iz}$, such that 
 $$
	F(z) = e^{zN}e^{\Gamma(s)}.F_{\infty}
$$
for $Im(z)>>0$.  
\endproclaim
\demo{Proof} One simply notes that the map $e^{-zN}.F(z)$ takes values in 
the homogeneous space $\check\M$, upon which $\Gc$ acts transitively.
To obtain $\Gamma(0) = 0$, note the definition of $F_{\infty}$.
\enddemo

\flushpar{\it Remark}. If $F(z)$ is unipotent in the sense of \cite{9}, 
i.e. the induced variations on $Gr^W$ are constant, the  
function $\Gamma(s)$ of $(2.3)$ assumes values in the 
in the subalgebra
$$
	Lie_{-1} = \{\, \a\in End(V) \mid \a:W_k\to W_{k-1}\}
	\subseteq Lie(\Gc)
$$ 

\par A structure of graded-polarized mixed Hodge 
structure $(F,W)$ on 
$V_{\R}\otimes\C$ induces one on the Lie 
algebra $Lie(\Gc)$.  Moreover, if $V = \bigoplus_{p,q}\, I^{p,q}$ is
the bigrading of $V$ defined by  $(F,W)$, then the corresponding decomposition of $Lie(\Gc)$ is given by 
the subspaces
$$
	gl(V)^{r,s} = \{ \a\in Lie(\Gc) \mid \a:I^{p,q}\to I^{r,s}\,\}
$$
Accordingly, $Lie(\Gc)=Lie(\Gc^F)\oplus q$, where 
$$ 
   Lie(G^F_{\Bbb C}) = \np\oplus\nz, \qquad
   q =\nn\oplus\lam\tag{2.4}
$$
and 
$$ \aligned
   \np &= \bigoplus_{r\geq 0,\,s<0}\,gl(V)^{r,s}, \\
   \nz &= gl(V)^{0,0},
   \endaligned\qquad
 \aligned
   &\nn = \bigoplus_{s\geq 0,\,r<0}\,gl(V)^{r,s}, \\
   &\lam = \bigoplus_{r,s<0}\,gl(V)^{r,s}.
 \endaligned
$$
The decomposition 
$$Lie(\Gc) = \np\oplus\nz\oplus\nn\oplus\lam$$ 
determines projection operators 
$$
\aligned
   \pip &:Lie(\Gc)\to\np, \\
   \piz &:Lie(\Gc)\to\nz,
\endaligned\qquad
\aligned
   \pin &:Lie(\Gc)\to\nn, \\
   \pil &:Lie(\Gc)\to\lam,
\endaligned\tag{2.5}
$$
The congruence $(1.1)$ is 
reflected in the following relations
$$ \alignedat 2
      &\bar\np\subseteq\nn\oplus\lam,\qquad 
      &&\bar\nz\subseteq\nz\oplus\lam,  \\ 
      &\bar\nn\subseteq\np\oplus\lam,\qquad 
      &&\overline{\lam}=\lam.
     \endalignedat\tag{2.6}   
$$

Note that the subalgebra $\lam$ depends only upon the mixed
Hodge structure $(F,W)$ and not the choice of graded-polarization.
\vskip 10pt

% Was 2.5
\proclaim{Lemma 2.7} (\cite{10}, \cite{14}) Let $(F,W)$ be a mixed Hodge 
structure.  Then, 
$$
	g\in\exp(\lam_{(F,W)}) \implies I^{p,q}_{(g.F,W)} = g.I^{p,q}_{(F,W)}
$$
In particular, $F\in\M$ and 
$ g\in G_{\R}\cup\exp(\lam_{(F,W)})$ implies $ I^{p,q}_{(g.F,W)} 
  = g.I^{p,q}_{(F,W)}$.
\endproclaim

\par Given $F\in\M$, we can smoothly decompose [10][14]
$$
	g_{\C} = g_{\R}e^{\l}f					\tag{2.8}
$$
for $g_{\C}\in \Gc$, $g\approx 1$, with $$g_{\R}\in G_{\R},\qquad e^{\l}\in\exp(\lam),\qquad f\in\Gc^F.$$
Similarly,
$$
	g = g_{\R}e^{\l}f				
$$
 for $g\in G$, with 
$$g_{\R}\in G_{\R},\qquad e^{\l}\in\exp(\lam), \qquad f\in G^F_{-1} = \Gc^F\cap\exp(Lie_{-1})$$ 
Therefore $(2.8)$
holds whenever $g_{\C}.F\in\M$.  

\par For $g\in\exp(Lie_{-1})$ we can be more precise. Define an {\it extended Lie monomial} 
to be a non-zero function $\mu:Lie_{-1}\to Lie_{-1}$ which can be obtained 
from the primitive monomial $m(\a) = \a$ in a finite number of steps via the 
operations of complex conjugation, taking Hodge components with respect to 
$(F,W)$ and forming Lie brackets 
(e.g. $\mu(\a) = [\bar\a,\a^{-1,0}]^{-1,-1}$).  In particular, if $\mu$ is 
an extended Lie monomial then so is any non-zero function of the form 
$\pip(\mu)$, $\pil(\mu)$ or $\pin(\mu).$ Define
$$
	Lie_{-r} = \{\, \zeta\in gl(V) \mid \zeta:W_k\to W_{k-r} \,\}
$$ 
which are ideals of $Lie_{-1}$.

% We also note that if $\a\in Lie_{-1}$ then
% $$
% 	\a = \pip(\a) + \pil(\a) + \pin(\a)
% $$

\proclaim{Theorem (2.9)} Let $F\in\M$  and $\a\in Lie_{-1}$.  Then,
$$
	e^{\a} = e^{\g}e^{\l}e^{\phi}
$$
where $\g$, $\l$ and $\phi$ are extended Lie polynomials in $\a$ which 
take values in the respective subalgebras 
$Lie(G_{\R})_{-1} = Lie(G_{\R})\cap Lie_{-1}$, $\lam$ and 
$Lie(\Gc^F)_{-1} = Lie(\Gc^F)\cap Lie_{-1}$. Moreover,
$$
\aligned
	\g   &= \pin(\a) + \overline{\pin(\a)} \mod Lie_{-2},		\\
	\l   &= \pil(\a) - \pil(\overline{\pin(\a)}) \mod Lie_{-2},	\\
	\phi &= \pip(\a) - \pip(\overline{\pin(\a)}) \mod Lie_{-2}.
\endaligned							\tag{2.10}
$$
\endproclaim
\demo{Proof} Since $Lie_{-1}$ is nilpotent, the 
Campbell--Baker--Hausdorff formula for it terminates 
after finitely many terms.  Set 
$$
	\g_1   = \pin(\a) + \overline{\pin(\a)},\qquad
	\l_1   = \pil(\a) - \pil(\overline{\pin(\a)}),\qquad
	\phi_1 = \pip(\a) -\pip(\overline{\pin(\a)}).
$$
Because of $(2.6)$ and
 $
	[Lie_{-p},Lie_{-q}]\subseteq Lie_{-p-q},
 $
one has $\g_1 + \l_1 + \phi_1 = \a$. Therefore
$$
	e^{\g_1}e^{\l_1}e^{\phi_1} = e^{\a + \b_1}
$$
where $\b_1$ is an extended Lie polynomial with values in $Lie_{-2}$. Inductively, suppose that there exist extended Lie polynomials
$\g_k$ with values in $Lie(G_{\R})_{-1}$, $\l_k$ with values in $\lam$
and $\phi_k$ with values in $Lie(\Gc^F)_{-1}$ such that
$$
	e^{\g_k}e^{\l_k}e^{\phi_k} = e^{\a + \b_k}
$$
for some extended Lie polynomial $\b_k$ with values in $Lie_{-k-1}$.  Let
$$
\aligned
	\g_{k+1}   &= \g_k + \g',\qquad 
	  \g' = -\pin(\b_k) - \overline{\pin(\b_k)}			\\
	\l_{k+1}   &= \l_k + \l',\qquad
	   \l'   = -\pil(\b_k) + \pil(\overline{\pin(\b_k)})		\\
	\phi_{k+1} &= \phi_k + \phi',\qquad  
	   \phi' = -\pip(\b_k) + \pip(\overline{\pin(\b_k)}
\endaligned
$$
Then $\g' + \l' + \phi' = -\b_k$, and hence
$$
	e^{\tilde\g_{k+1}}e^{\l_{k+1}}e^{\phi_{k+1}}
	= e^{\a + \g' + \l' + \phi' + \b + \b_{k+1}}
	= e^{\a + \b_{k+1}}
$$
for some extended Lie polynomial $\b_{k+1}\in Lie_{-k-2}$.  
As $Lie_{-r} = 0$ for some index $r=r_0$, this completes the proof.
\enddemo
 
%%%%%%%%%%%%%%%%%%%%%%%%%%%% End section-2 %%%%%%%%%%%%%%%%%%%%%%%%%%%%%%

%%%%%%%%%%%%%%%%%%%%%%%%%%%%%% section-3 %%%%%%%%%%%%%%%%%%%%%%%%%%%%%%%%
$$
	\text{\bf 3.\quad Unipotent Variations}
$$

\par Let $X$ let be a smooth, complex algebraic variety.  Then, by the
work of Hain, Morgan et\. al\. (cf\. \cite{1} for an overview) for each 
positive integer $k$, there exists a canonical admissible, variation of 
graded-polarized mixed Hodge structure $\V\to X$ with fiber 
$$
	\V_x = \C\pi_1(X,x)/(J_x)^{k+1}
$$
where $J_x$ is the augmentation ideal of $\C\pi_1(X,x)$.  The monodromy representation
$$
	\rho:\pi_1(X,x)\to Aut(\V_x)
$$
of such a variation is unipotent and the variations of pure, polarized Hodge structure 
induced by $\V$ on $Gr^{\W}$ are constant in this situation. Such variations are called {\it unipotent}. 
With this motivation in mind, we consider here
problem of determining the admissibility for such variations
in terms of the grading $\Y$ as discussed in [\S 1], begining with the following observations:
\roster
\item"(a)" By virtue of Schmid's Nilpotent Orbit Theorem, the
monodromy logarithm $N$ of a unipotent variation $\V\to\Delta^*$ must act
trivially on $Gr^W$.
\item"(b)" On account of $(a)$, the relative weight filtration $\rel W$ of 
a unipotent variation exists if and only if $N:W_k\to W_{k-2}$ for each 
index $k$, i.e. $\rel W = W$.
\endroster

\proclaim{Theorem 3.1} Let $\V\to\Delta^*$ be a unipotent variation of 
graded-polarized mixed Hodge structure which is admissible.  Then, $\Y$ extends continuously, as a grading of $\W$, to 
the canonical extension $\tilde\V\to\Delta$.
\endproclaim
\demo{Proof} Let $F:U\to\M$ be a lift the period map of $\V$ to the 
upper half-plane $U$ and let $N$ be its monodromy logarithm.  Then, as discussed in \S 2, the claim is 
equivalent to the assertion that 
$$
	Y_{\infty} = \lim_{\text{Im}(z)\to\infty}\, e^{-zN}.Y_{(F(z),W)}
	\tag{3.2}
$$
exists and grades $W$.  

\par To verify the existence of the grading $Y_{\infty}$, observe that on 
account of the unipotency of $\V$, the function 
$$
	\psi(s) = e^{-zN}.F(z),\qquad s = e^{2\pi iz}	
$$
considered in \S 2 takes values in $\M$ (and not $\check\M$ as is a priori
the case).  Indeed, in our case $\psi(s)$ induces the same filtration as $F(z)$ on $Gr^W$, and 
hence is an element of $\M(Gr^W)$.  Moreover, since the induced filtration 
$F(z)Gr^W$ of a unipotent variation is by definition constant, the limiting
Hodge filtration
$$
	F_{\infty} = \lim_{s\to 0}\,\psi(s)
$$
of a $\V$ is likewise an element of $\M$.  
By $(b)$ above, the relative weight
filtration of our admissible variation $\V$ must coincide with the weight
filtration $W$ of $\V$.  By Deligne's theorem,  $N$ is $(-1,-1)$--morphism of 
$(F_{\infty},W)$.
We now apply Lemma $(2.3)$, with $q = \nn\oplus\lam$ as
defined by  $(2.4)$ and $F = F_{\infty}$.  This gives
$$
	F(z) = e^{zN}e^{\Gamma(s)}.F_{\infty},\qquad \Gamma(0) = 0									
$$
relative to a holomorphic function $\Gamma(s)$ taking values in the subalgebra
$$
	q\cap Lie_{-1} = \bigoplus_{r<0,r+s\leq -1}\, gl(V)^{r,s}
								\tag{3.3}
$$
Consequently,
$$
	e^{-zN}.Y_{(F(z),W)} = e^{-zN}.Y_{(e^{zN}e^{\Gamma(s)}.F_{\infty},W)}
$$

\par To finish the proof, let $e^{X(z)} = e^{zN}e^{\Gamma(s)}$  Then 
$$
	e^{X(z)} = e^{zN + \Gamma(s) + 
		      (\text{brackets of $zN$ and $\Gamma(s)$})}
								\tag{3.4}
$$
while, by Theorem $(2.9)$,
$$
	e^{X(z)} = e^{\g(z)}e^{\l(z)}e^{\phi(z)}
$$
where $\g(z)$, $\l(z)$ and $\phi(z)$ extended Lie polynomials in $X(z)$
with respect to $(F_{\infty},W)$.
Accordingly, 
$$
\aligned
	e^{-zN}.Y_{(e^{zN}e^{\Gamma(s)}.F_{\infty},W)}
	&= e^{-zN}e^{\g(z)}e^{\l(z)}.Y_{(F_{\infty},W)}	\\
	&= e^{-zN}e^{zN}e^{\Gamma(s)}e^{-\phi(z)}.Y_{(F_{\infty},W)}	\\
	&= e^{\Gamma(s)}e^{-\phi(z)}.Y_{(F_{\infty},W)}
\endaligned							\tag{3.5}
$$
Since $\phi(z)$ is an extended Lie polynomial in $X(z)$ which takes
values in $Lie(\Gc^{F_{\infty}})_{-1}$, the Hodge components $X(z)^{r,s}$ 
of $X(z)$ with both $r$, $s<0$ can only appear in $\phi(z)$ inside of Lie
brackets which contain some Hodge component $X(z)^{p,q}$ with either 
$p$ or $q$ greater than zero.  Indeed, since
$$
	\lam = \bigoplus_{r,s<0}\, gl(V)^{r,s}
$$
is a subalgebra of $Lie_{-1}$ which is closed under complex conjugation, 
and the extended Lie monomials are exactly the non-zero functions 
$\mu:Lie_{-1}\to Lie_{-1}$ which can by constructed in a finite number
of steps from the primitive monomial $m(\a) = \a$ via the operations of
complex conjugation, taking Hodge components and forming Lie brackets,
any Lie monomial $\mu(\a)$ which depends only on the Hodge components 
$\a^{r,s}$, with both $r$, $s<0$ must preserve $\lam$.  

\par Returning to  $(3.4)$ and noting that $N\in gl(V)^{-1,-1}$, it follows from the above 
considerations that 
$$
	||\phi(z)|| \leq K|z|^b |e^{-2\pi iz}|
$$
for some constants $K$ and $b$ (and $||\ast||$ a norm on $Lie_{-1}$).
Consequently, 
$$
	\lim_{Im(z)\to\infty}\, \phi(z) = 0.			\tag{3.6}
$$
On
account of  $(3.5)$, $(3.6)$ and the fact that 
$\Gamma(0) = 0$, the limit $(3.2)$ exists and equals $Y_{(F_{\infty},W)}$. 
\enddemo

This finishes the proof of $(3.1)$.

\proclaim{Corollary 3.7} Let $\V\to\Delta^*$ be a unipotent variation of 
graded-polarized mixed Hodge structure which is admissible.  Then, $\Y$ extends continuously, as a grading of $\W$, to 
the canonical prolongation of $\V$, with limiting value 
$Y_{\infty} = Y_{(F_{\infty},W)}$.
\endproclaim

\par We now establish the converse of Theorem $(3.1)$:

\proclaim{Theorem (3.8)} Let $\V\to\Delta^*$ be a unipotent variation
of graded-polarized mixed Hodge structure, and suppose that the grading
$\Y$ of the underlying weight filtration $\W$ defined by the $I^{p,q}$'s
of $\V$ extends continuously, as a grading of $\W$, to the canonical 
extension $\tilde\V$.  Then, $\V$ is admissible.
\endproclaim
\demo{Proof} As in the proof of Theorem $(3.1)$, we select a 
lifting $F:U\to\M$ of the period map of $\V$.  Lemma $(2.1)$,
implies that $\Y$ extends continuously to a
grading of $\W$ in the canonical extension of $V$, that
$$
	Y_{\infty} = \lim_{s\to 0}\, Y(s),\qquad
	Y(s) = e^{-zN}.Y_{(F(z),W)}			\tag{3.9}
$$
with $s=e^{2\pi iz}$, exists and grades $W$. 

\par To see that $(3.9)$ implies the existence of the limiting Hodge
filtration 
$$
	F_{\infty} = \lim_{s\to 0}\, \psi(s),\qquad \psi(s) = e^{-zN}.F(z)
$$
observe that  $Y_{(F(z),W)}$ preserves the filtration $F(z)$, so  $Y(s)$  preserves $\psi(s)$.
Accordingly, $\psi(s)$ can be obtained by pulling back the induced 
filtration $\psi(s) Gr^W$ via the isomorphism $V\cong Gr^W$ 
determined by $Y(s)$, namely
$$
	Gr^W_k \cong E_k(Y(s))
$$
As the induced filtrations $\psi(s)Gr^W$ are constant,  the existence of the limiting grading 
$(3.9)$ imply the existence of the limiting Hodge filtration 
$F_{\infty}$.

\par By (b), the existence of the relative weight filtration is equivalent (in the present case) to the statement $$N:W_k\to W_{k-2}.$$
As in the proof of Theorem $(3.1)$, the unipotency of $\V$ implies
that $F_{\infty}\in \M$. We can then apply Theorem $(2.3)$
with $F=F_{\infty}$ and $q = \nn\oplus\lam$.  This gives 
$$
	F(z) = e^{zN}e^{\Gamma(s)}.F_{\infty},\qquad \Gamma(s) = 0
$$
relative to a holomorphic function $\Gamma(s)$ taking values in the
subalgebra $(3.3)$.
Letting $$e^{X(z)} = e^{zN}e^{\Gamma(s)},$$ 
we may write (cf\.  $(3.5)$) 
$$
	Y(s) = e^{\Gamma(s)}e^{-\phi(z)}.Y_{(F_{\infty},W)}
$$ 
where $\phi(z)$ is the extended Lie polynomial in $X(z)$ obtained by 
decomposing $e^{X(z)}$ in accord with $(2.9)$.  Moreover,
by  $(2.10)$,
$$
	\phi(z) = \pip(X(z)) - \pip(\overline{\pin(X(z))}) \mod Lie_{-2}
							   \tag{3.10}
$$

\par Claim:
$$
	N:F_{\infty}^p\to F^{p-1}_{\infty}			\tag{3.11}
$$
Indeed, by the horizontality of $F(z)$,
$$
	\frac{d}{dz} F^p(z) \subseteq F^{p-1}(z)		
$$
Inserting the formula $F(z) = e^{zN}e^{\Gamma(s)}$ in the last equality and symplifying, one then obtains
$$
	e^{-\Gamma(s)}.N 
	  + 2\pi i s e^{-\Gamma(s)}\frac{d}{ds} e^{-\Gamma(s)}
	:F_{\infty}^p \to F_{\infty}^{p-1}			
$$
Setting $s=0$,  $(3.11)$ follows.

\par Because of  $(3.11)$, the mixed-Hodge decomposition of 
$N$ relative to $(F_{\infty},W)$ satifies
$$
	N = N^{0,-1} + N^{-1,0} \mod Lie_{-2}			\tag{3.12}
$$
with
$$
	\overline{N^{-1,0}} = N^{0,-1} \mod Lie_{-2}		\tag{3.13}
$$
In addition, by $(3.4)$,
$$
	X(z) = zN + \Gamma(s) \mod Lie_{-2}			\tag{3.14}
$$

\par Taking note of $(3.10)$ and $(3.12)$--$(3.14)$, it then follows that
$$
	\phi(z) = 2iyN^{0,-1} 
		  + \pip(\Gamma(s)) - \pip(\overline{\pin(\Gamma(s))})
		  \mod Lie_{-2}		
								\tag{3.15}
$$
Accordingly, the limiting grading 
$$
\aligned
	Y_{\infty} &= \lim_{s\to 0}\, Y(s)				   
		    = \lim_{\text{Im}(z)\to\infty}\, e^{-zN}.Y_{(F(z),W)}  \\
		   &= \lim_{\text{Im}(z)\to\infty}\, 
			e^{\Gamma(s)}e^{-\phi(z)}.Y_{\infty}
\endaligned							\tag{3.16}
$$
exists only if $N^{0,-1} = 0$.  Indeed, by \cite{3}, the group 
$\exp(Lie_{-1})$ acts simply transitively upon the set of all gradings $Y$ of 
$W$.  Therefore, in order for the limit $(3.16)$ to exist, $||\phi(z)||$ must
remain bounded as $Im(z)\to\infty$, and hence $N^{0,-1}$ must be equal to 
zero by  $(3.16)$ [recall $\Gamma(0) = 0$].  Since 
$N^{-1,0} = \overline{N^{0,-1}}$ by  $(3.13)$, it then follows
that $N = 0 \mod Lie_{-2}$, i.e. $N:W_k\to W_{k-2}$.
\enddemo

\par Combining Theorems $(3.1)$ and $(3.8)$, we then obtain:

\proclaim{Theorem (I)} Let $\V\to\Delta^*$ be a variation of graded-polarized
mixed Hodge structure which is unipotent.  Then, $\V$ is admissible if
and only if the $C^{\infty}$-grading $\Y$ of $\W$ defined by the $I^{p,q}$'s
of $\V$ extends continuously, as a grading of $\W$, to the canonical 
extension $\tilde\V$.
\endproclaim

%%%%%%%%%%%%%%%%%%%%%%%%%%%% End section-3  %%%%%%%%%%%%%%%%%%%%%%%%%%%%%%

%%%%%%%%%%%%%%%%%%%%%%%%%%%%%%% section-4 %%%%%%%%%%%%%%%%%%%%%%%%%%%%%%%%

$$
	\text{\bf 4.\quad Split Orbits}
$$

\par Let $\M = \M(W,\cal S, h^{p,q})$ be a classifying space of 
graded-polarized mixed Hodge structure, with \lq\lq compact dual\rq\rq{}
$\check\M$ and associated Lie groups $G_{\R}$, $G$ and $\Gc$, as described in section 2.

\definition{Definition} A nilpotent orbit of graded-polarized mixed
Hodge structure (modeled on $\M$) consists of a filtration $F\in\check\M$ 
and an element $N\in Lie(G_{\R})$ such that 
\roster
\itemb $N:F^p\to F^{p-1}$ for each index $p$.
\itemb There exists a constant $\a$ such that 
       $\text{Im}(z)>\a \implies e^{zN}.F\in\M$.
\endroster
\enddefinition

\par In this section  we prove the following version of Theorem II:

\proclaim{Theorem 4.1} A nilpotent orbit of graded-polarized mixed Hodge
structure $e^{zN}.F$ is admissible and split if and only if
\roster
\item"(a)" The limit
$$
 	Y_{\infty} = \lim_{Im(z)\to\infty}\, Ad(e^{-zN}) Y_{(e^{zN}.F,W)}
$$
exists, grades $W$ and is defined over $\R$.
\item"(b)" Each of the induced orbits $e^{zN}.F Gr^W_k$ is split
\endroster
\endproclaim

\par As noted in the introduction, the fact that an admissible split nilpotent orbit satisfies the conditions  of the Theorem  is an 
immediate consequence of results of \cite{6} discussed in the attached 
appendix.

\par To prove the converse, we let 
$$
	\rho:sl_2(\R) \to End_{\R}(V)
$$
denote the representation of $sl_2(\R)$ constructed in \S 1 by simply 
pulling back the representations $\rho_k:sl_2(\R)\to End_{\R}(Gr^W)$
to $V$ via the isomorphism $Gr^W\cong V$ induced by $Y_{\infty}$, and
define
$$
	\prel Y = Y_{\infty} + Y_o,\qquad Y_o = \rho(h)
$$
  Then, as discussed in \S 1, in order
for Theorem (4.1) to hold, $\prel Y$ must grade the relative weight
filtration of $N$ and $W$, i.e. the associated filtration
$$
	\prel W_k = \bigoplus_{\ell\leq k}\, E_j(\prel Y)
$$
must satisfy the following two conditions (cf. $(2.2)$):
\roster
\item"(i)"  $N:\prel W_k\to\prel W_{k-2}$.
\item"(ii)" $\prel W$ induces on $Gr^W_k$ the shifted monodromy weight 
filtration
$$
	W^{\#}(k) = W(N:Gr^W_k\to Gr^W)[-k]
$$
of $N$ on $Gr^W$.
\endroster

\par To  facilitate the proof we  record the
following observations: 
\roster
\item The limiting grading $Y_{\infty}$  from (a) in Theorem (4.1), preserves $F$.

\item $Y_o=\rho(h)$ preserves $F$. 

\item Suppose that $\prel Y$ does indeed grade the relative weight filtration
of $N$ and $W$.  Then $(F,\rel W)$ is split over $\R$.  In particular, the
proof of Theorem (4.1) will be complete once we verify  (i) and 
(ii) above.

\item The representation $\rho$ defined above acts on $Gr^W$ by infinitesimal
isometries (i.e\. $\rho$ takes values in $Lie(G_{\R})$.
\endroster
To prove (1), one simply notes that $e^{-zN}.Y_{(e^{zN}.F,W)}$ preserves $F$ 
whenever $(e^{zN}.F,W)$ is a mixed Hodge structure.  To prove (2), we simply 
recall from \cite{15} that:
$$
	\rho_k(h) = Y_{(FGr^W_k,W^{\#}(k))} - k
$$
To prove (3), observe that items (1) and (2) together with the hypothesis of Theorem, imply that $\prel Y = Y_{\infty} + Y_o$ is a semisimple
endomorphism of $V$ which preserves $F$.  Accordingly, if $\prel Y$ also 
grades $\rel W$ then $(F,\rel W)$ must be a split mixed Hodge structure.  The proof of (4) can be found in \cite{12} and boils down to the assertion that the monodromy weight filtration 
$W(N:Gr^W_k\to Gr^W_k)$ is self-dual with respect to the polarization ${\S}_k$.
\vskip 10pt

\par To verify  (ii), one simply observes that on account of the 
formula for $\rho_k$ given above, $\prel Y = Y_{\infty} + \rho(h)$ acts
on $Gr^W_k$ as $Y_{(F,W^{\#}(k))}$.

\par To establish  $(i)$, we shall actually verify a (seemingly) 
stronger condition, namely  
$$
	[\prel Y,N] = -2N			\tag{4.2} % was 5.1
$$ 

\par To this end, we let
$$
	N = N_0 + N_{-1} + \cdots 
$$
denote the decomposition of $N$ with respect to the eigenvalues of 
$ad\,Y_{\infty}$ and recall that in the present context:
$$
	Y_o = \rho(h),\qquad N_0 = \rho(n_-)
$$
Accordingly,  $(4.2)$ holds if and only if 
$$
	[Y_o,N_{-k}] = (k-2)N_{-k}		
$$
for each index $k>0$.  Consequently, it will suffice to prove the following
assertion:

% Lemma (4.3) was (5.3)
\proclaim{Lemma 4.2} Relative to the representation $\rho$ of $sl_2(\R)$ 
defined above, any non-zero component $N_{-k}$ with $k>0$ is a highest 
weight vector of weight exactly $(k-2)$.
\endproclaim

For the remainder of this section, we
shall write $Y_F$ in place of $Y_{(F,W)}$ whenever $F\in\M$.
\vskip 10pt

\par By virtue of  (a), we know that:
$$
	Y_{\infty} = \lim_{y\to\infty}\,e^{-iyN}.Y_{e^{iyN}.F}
							\tag{4.3} % was 5.4
$$
Furthermore, we also know that $N$ is horizontal with respect to $F$, 
and hence
$$
	N_{-k}:F^p \to F^{p-1}				\tag{4.4} 
$$
since $Y_{\infty}$ preserves $F$.  Accordingly, we can  prove
Lemma $(4.2)$ inductively by explicitly computing the right hand side of 
 $(4.3)$ and then imposing the horizontality condition $(4.4)$.
\vskip 10pt

\par To fill in the details, observe that  Schmid's $SL_2$ Orbit Theorem applied to $Gr^W$ (a sum of variations of pure structures) together with  (b), implies that 
the pair
$$
	(e^{iyN_0}.F,W)
$$
is a mixed Hodge structure for all $y>0$.  In particular, the base point
$$
	F_o = e^{iN_0}.F		 % was 5.6
$$
is an element of $\M$.  Moreover, since $Y_{\infty}$ is both real and
preserves the filtration $e^{iyN_0}.F$,
$$
	Y_{\infty} = Y_{(e^{iyN_0}.F,W)}
$$
for all $y>0$.

\par Next, we introduce the $\exp(Lie_{-1})$-valued function
$$
	e^{Q(y)} = e^{iyN}e^{-iyN_0} = e^{iyN_0 + iyN'}e^{-iyN_0},\qquad
	N' = N - N_0  					\tag{4.5} % was 5.7
$$
and note that
$$
	e^{iyN} = e^{Q(y)}e^{iyN_0} 
		= y^{-\half Y_o}e^{P(y)}e^{iN_0}y^{\half Y_o}	\tag{4.7}
$$
upon setting $P(y) = Ad(y^{\half Y_o})Q(y)$.  In particular, since
$y^{\half Y_o}$ preserves $F$:
$$
	e^{iyN}.F 
	= y^{-\half Y_o}e^{P(y)}e^{iN_0}y^{\half Y_o}.F
	=  y^{-\half Y_o}e^{P(y)}.F_o
$$
Thus, upon applying Theorem $(2.9)$ to $e^{P(y)}$ with $F=F_o$, we obtain
a distinguished decomposition
$$
	e^{P(y)} = e^{\g(y)}e^{\l(y)}f(y),\qquad 
	f(y) = e^{\phi(y)} 	
							 % was 5.8
$$
such that
$$
	I^{p,q}_{(e^{P(y)}.F_o,W)} = e^{P(y)}f^{-1}(y).I^{p,q}_{(F_o,W)}
$$
Consequently,
$$
\aligned
	Y_{e^{iyN}.F} 
	&= y^{-\half Y_o}.Y_{e^{P(y)}.F_o}			
	 = y^{-\half Y_o}e^{P(y)}f^{-1}(y).Y_{(F_o,W)}		\\
	&= y^{-\half Y_o}e^{P(y)}f^{-1}(y).Y_{\infty}
\endaligned
$$
wherefrom [cf\. (4.7)]
$$
\aligned
	e^{-iyN}.Y_{e^{iyN}.F} 
	&= e^{-iyN}y^{-\half Y_o}e^{P(y)}f^{-1}(y).Y_{\infty} 	\\
	&= (y^{-\half Y_o}e^{P(y)}e^{iN_0}y^{\half Y_o})^{-1}
	   y^{-\half Y_o}e^{P(y)}f^{-1}(y).Y_{\infty}		\\
	&= y^{-\half Y_o}e^{-iN_0}f^{-1}.Y_{\infty}		\\
	&= y^{-\half Y_o}e^{-iN_0}f^{-1}(y)e^{iN_0}.Y_{\infty}  
\endaligned
$$
with the very last step being justified by the fact that 
$e^{iN_0}.Y_{\infty} = Y_{\infty}$.

\par Since $[Y_o,Y_{\infty}]=0$, 
$$
	e^{-iyN}.Y_{e^{iyN}.F} 
	= [Ad(y^{-\half Y_o})Ad(e^{-iN_o})f^{-1}(y)].Y_{\infty}
$$
In particular, since $f^{-1}(y)$ takes values in $\exp(Lie_{-1})$, 
$(4.3)$ holds if and only if 
$$
	\lim_{y\to\infty}\, Ad(y^{-\half Y_o})Ad(e^{-iN_o})f^{-1}(y) = 1
$$
Equivalently, if $f(y) = e^{B(y)}$, then
$$
	\lim_{y\to\infty}\, Ad(y^{-\half Y_o})Ad(e^{-iN_o})B(y) = 0
							\tag{4.7} 
$$

\par Write 
$$
	B(y) = \sum\, B_m y^{m\over 2} 
$$
Then, 
$$
	Ad(y^{-\half Y_o})Ad(e^{-iN_0})B(y) 
	= \sum_m\sum_j [Ad(e^{-iN_0})B_m]^{Y_o}_jy^{m-j\over 2}
$$
and hence  $(4.7)$ holds if and only if
$$
	[e^{-i ad\,N_0}\,B_m]^{Y_o}_j = 0, 
	\qquad \forall\hph{a}j\leq m 			\tag{4.8}
$$
Let $[A]^S_{\lambda}$  denote the component of $A$ in the 
$\lambda$-eigenspace of a semisimple element $S\in End(V)$.
\vskip 10pt

\par To compute $B(y)$, note that \cite{6} (cf. Appendix):  
$$
\aligned 
       e^{Q(y)} &=e^{(iyN_0+iyN')}e^{-iyN_0} \\
 		&=\exp(\Psi(iyN', \dots, {1\over m!}ad(iyN_0)^m(iyN'),\dots))
\endaligned
$$
for some  Lie polynomial $\Psi$ such that
$$ 
	Q(y) = {e^{ad(iyN_0)}-1\over ad(iyN_0)}(iyN') \ + \ O((iyN')^2)	
$$
where $O((iyN')^2)$ is of ``lower weight", i.e. if $N'=N_{-k}$,  this is
a decomposition of $Q(y)$ according to $E_{-k}(Y_{\infty}) + Lie_{-2k}$.
Consequently,
$$ 
	P(y)=   \frac{e^{ad(iN_0)}-1}{ad(iN_0)}(iN'(y)) \ 
		+ \  U_{[-2]}(y) 			    \tag{4.9}
$$
where 
$$ 
   N'(y)= y Ad(y^{Y_o\over 2}) N' 
	= \sum_{k\geq 1}\sum_j\, [N_{-k}]^{Y_o}_j\, y^{1+\half j}
							    \tag{4.10}
$$ 
and $U_{[-2]} = y^2 Ad(y^{\frac{Y_o}{2}})O((N')^2)$ is of lower weight.
\vskip 10pt

\par Now, $Lie(\Gc)$ decomposes
as a direct sum of irreducible $sl_2$-modules $U$ under the adjoint action
of $\rho$.  Moreover, since $Y_{\infty}$ is both defined over $\R$ and 
commutes with $\rho$, $U$ is likewise defined over $\R$ and contained in 
some eigenspace of $ad\,Y_{\infty}$.  
In order to better understand that decomposition, let us assume 
$U\subseteq Lie(\Gc)$ to be an irreducible submodule contained in the 
$-k$ eigenspace of $ad\,Y_{\infty}$ for some positive value of $k$.  
To see that $U$ inherits a pure Hodge structure of weight $-k$ from 
$$
	Lie(\Gc) = \bigoplus_{r,s}\, gl(V)^{r,s}_{(F_o,W)}	\tag{4.11}
$$
we observe that in the present context, the adjoint 
representation $\rho:sl_2(\C)\to Lie(\Gc)$ becomes a morphism of mixed 
Hodge structures, upon endowing $sl_2(\C)$ with the following pure Hodge 
structure of weight zero:
$$
	x^- = \half(h - in_- - in_+)\in H^{-1,1},\qquad
	x^+ = \half(h + in_- + in_+)\in H^{1,-1},
$$
$$
	z   = i(n_- - n_+)\in H^{0,0},
$$	
with $N_o = \rho(n_-)$ and 
$Y_o = \rho(h)$.  Since the action of $Z = \rho(z)$ stabilizes 
each summand of $(4.17)$, the Hodge decomposition 
$$
	\cdots + gl(V)^{0,-k} 
	       + (gl(V)^{-1,-k+1} + \cdots + gl(V)^{-k+1,-1}) 
	       + gl(V)^{-k,0} + \cdots
	\tag{4.12}
$$
of $U$ must be exactly parallel to the corresponding decomposition
$$
	\cdots + E^Z_k + (E^Z_{k-2} + \cdots + E^Z_{-k+2}) + E^Z_{-k} + \cdots
	\tag{4.13}
$$
of $U$ with respect to the eigenvalues of $\rho(Z)$.
\vskip 10pt

\flushpar{\it Remark.} To avoid confusion regarding the Hodge
types of the summands appearing in the preceding equation, we note that 
$(4.10)$ is equivalent to the assertion that $E^Z_{k+2p}$ has Hodge type 
$(p,-k-p)$.  We also note that for $k>1$, the parenthesized terms appearing
in  $(4.10)$ correspond to $U\cap\lam_{(F_0,W)}$.
\vskip 10pt

\par We may now proceed with
the proof of Lemma $(4.2)$.  Inductively, we may assume that it holds for eigenvalues of $ad\, Y_{\infty}$ bigger than $-k$.  In 
particular, since the $su_2$-basis $(X^-,Z,X^+)$ given by 
$$
	X^- = \rho(x^-),\qquad  Z = \rho(z),\qquad  X^+ = \rho(x^+)
$$ 
defines the same representation of $sl_2(\C)$ as $(N_0,Y_o,N^+_0)$, the
 monomials 
$$
	(ad\,N_0)^j\,N_{-\ell},\qquad -\ell > -k, \quad 0\leq j \leq \ell-2
								\tag{4.14}
$$ 
belong to $\lam_{(F_0,W)}$ on account of  $(4.10)$.  Consequently, if
we decompose 
$$
	P(y) = P_{-1}(y) + P_{-2}(y) + \cdots
$$
according to the eigenvalues of $ad\, Y_{\infty}$, $(4.14)$ and  $(4.11)$ imply 
that 
$$
	P(y) = P_{\R}(y)+ P_{\Lambda}(y) + B(y) \mod Lie_{-k-1}  \tag{4.15}
$$
with $P_{\R}(y)\in Lie(G_{\R})_{-1}$, $P_{\Lambda}(y)\in\lam_{(F_o,W)}$ and
$$
	P_{-k}(y) 
	 =  \frac{e^{ad(iN_0)}-1}{ad(iN_0)}
	     (i\sum_m\, y^{1+\half m}  [N_{-k}]^{Y_o}_m) \mod \lam_{(F_o,W)}
								\tag{4.16}
$$
Indeed: $(4.11)$ gives a formula for $P(y)$ in terms of the monomials
$(ad\,N_0)^r N_{-s}$ and their projections.  By inductive 
hypothesis and the fact that $\lam_{(F_o,W)}$ is an associative subalgebra
of $gl(V)$, the monomials listed in $(4.14)$ can only contribute factors 
which belong to $\lam_{(F_o,W)}$.

\par Together with $(4.12)$, this shows
$$
	B(y) = \text{proj}_{F_o}\left(
		  \frac{e^{ad(iN_0)}-1}{ad(iN_0)}
	            \sum_m\, i y^{1+\half m}  [N_{-k}]^{Y_o}_m\right)
		\mod Lie_{-k-1}					\tag{4.17}
$$
where, in general, we shall let 
$$
	\text{proj}_{F_o}: Lie(\Gc) \to Lie(\Gc^{F_o})
$$
denote the projection operator defined by  $(4.11)$ via the rule:
$$
	\text{proj}_{F_o}(\sum_{r,s}\, \a^{r,s}) 
	= \sum_{r\geq 0, s}\, \a^{r,s}
$$ 
Accordingly, by  $(4.17)$, the coefficients $B_m$ of the series 
expansion  $B(y) = \sum_m\, B_m y^{\half m}$ are given by 
$$
   B_m = proj_{F_o}\Big(\frac{e^{ad(iN_0)}-1}{ad(N_0)}[N_{-k}]^{Y_o}_{m-2}\Big)
	 \mod Lie_{-1-k}					
$$

\par Returning now to the setting of $(4.10)$, let us restrict our attention
to an irreducible $sl_2$-module 
$$
	U\subseteq E_{-k}(ad\,Y_{\infty})
$$
of highest weight $d$.  Then, by virtue of our preceding remarks, we obtain
the identity
$$
  B_m = proj_{F_o}\Big(\frac{e^{ad(iN_0)}-1}{ad(N_0)}[N_{-k}]^{Y_o}_{m-2}\Big) $$

\par To continue, we note that for positive values of $k$, we may compute
$proj_{F_o}(T)$ for any element $T\in E_{-k}(Y_{\infty})$ via the formula
$$
   Lie(\Gc) = Lie(G_{\R}) + \lam_{(F_o,W)} + Lie(\Gc^{F_o})
$$
by determining its Hodge components: replacing $T$ by $\half(T-\bar T)$ gets 
rid of the real part, and taking just the components of type 
$..., (1,-k-1), (0,-k)$ removes the contribution from $\Lambda^{-1,-1}$.  
Thus,
$$
	proj_{F_o}(T) = \half \sum_{a\geq 0} (T  - \bar T)^{a,-k-a}
$$
In particular,
$$
\aligned
  B_m 
   &= proj_{F_o}\Big(\frac{e^{ad(iN_0)}-1}{ad(N_0)}[N_{-k}]^{Y_o}_{m-2}\Big) \\
   &= \half\sum_{a\geq 0}\, 
      \Big((\frac{e^{ad(iN_0)}-1}{ad(N_0)} 
	    - \frac{e^{-ad(iN_0)}-1}{ad(N_0)})[N_{-k}]^{Y_o}_{m-2}
      \Big)^{a,-k-a}							     \\
   &=  i\sum_{a\geq 0}\, 
         \Big(\frac{\sin ad(N_0)}{ad(N_0)}[N_{-k}]^{Y_o}_{m-2}
         \Big)^{a,-k-a}
\endaligned
$$
Therefore, by  $(4.10)$:
$$
	B_m = i\sum_{\ell\geq k}\, 
         \Big[\frac{\sin ad(N_0)}{ad(N_0)}[N_{-k}]^{Y_o}_{m-2}
         \Big]^Z_{\ell}					\tag{4.18}
$$

\par Since
$F_o = e^{iN_0}.F$, the horizontality of $N_{-k}$ at $F$ together with $(4.10)$
imply that 
$$
	e^{i ad\, N_0}N_{-k}
	\in\bigoplus_{a\geq -1}\, gl(V)^{a,-k-a}_{(F_o,W)}
	  =\bigoplus_{\ell\geq k-2}\, E^Z_{\ell}	\tag{4.19}
$$
and hence 
$$d\geq k-2.$$
\par 
 $(4.19)$ is equivalent to the assertion that 
$$
	N_{-k}\in\bigoplus_{\ell\geq k-2}\, E^{Y_o}_{\ell},	\tag{4.20}
$$
since $e^{i ad\, N_0}$ is actually an isomorphism from $E^{Y_o}_{\ell}$ to
$E^Z_{\ell}$ (see proof of Lemma $(4.23)$ below).

\proclaim{Corollary} By virtue of horizontality, $B_m = 0$ for $m<k$.
\endproclaim
\demo{Proof} Indeed, by  $(4.20)$, $[N_{-k}]^{Y_o}_j = 0$ unless
$j\geq k-2$.  Therefore, by  $(4.18)$, $B_m = 0$ unless $m\geq k$.
\enddemo

\par To prove that the remaining coefficients $B_m = 0$, we note that by 
 $(4.8)$, in order for the limiting grading $Y_{\infty}$ to exist,
we must have
$$
	e^{-i ad\, N_0} B_j \in \bigoplus_{\ell\geq j+2}\, E^{Y_o}_{\ell}
$$
for any non-zero coefficient $B_j$, and hence 
$$
	B_j \in \bigoplus_{\ell\geq j+2}\, E^Z_{\ell}	       \tag{4.21}
$$
since $e^{i ad\, N_0}$ is an isomorphism from $E^{Y_o}_{\ell}$ to $E^Z_{\ell}$.

\proclaim{Corollary} $B_j = 0$ if $j\geq d$.
\endproclaim

\par To dispense with the remaining coefficients 
$$
	B_k, \dots, B_{d-2}					\tag{4.22}
$$
we will use that,  
in the present context, 
$$
	\Big[{\sin{N_0} \over N_0 } E^{Y_o}_a\Big]^Z_{b}\not=0
$$
iff both $a$ and $b$ are weights and  either $a$ is a highest weight, 
or  $a<|b|$. This will be proved in $(4.23)$ below.
\vskip 10pt

\par To see that all of the coefficients  in $(4.22)$ must then 
vanish, note that by  $(4.18)$, such a coefficient $B_j$ could
be non-zero if and only if:
$$
	[N_{-k}]^{Y_o}_{j-2} \neq 0			        
$$
If $d$ was bigger than $d-2$, $(4.23)$ would imply
$$
	[B_j]^Z_j \neq 0
$$
on account of  $(4.18)$.  As this contradicts $(4.21)$, all the
coefficients listed in $(4.22)$ must vanish.  Accordingly, $N_{-k}$ must
be of exactly highest weight $k-2$. 
\vskip 10pt

It only remains to prove

% Lemma 4.23, was 4A.
\proclaim{Lemma 4.23}
Let $(\rho,V)$ be a finite-dimensional representation of $sl_2({C})$,
$$
Y   = \rho \pmatrix 
  	1 & 0  \\
  	0 & -1 \\ 
        \endpmatrix,\qquad Z   = 
\rho \pmatrix
  	0 & -i \\
  	i & 0 \\
         \endpmatrix, \qquad 
N   = \rho\pmatrix 
  	0 & 0 \\
 	1 & 0 \\
        \endpmatrix,$$
$y_a$, $z_a$, the natural projections to the $a$-weight spaces of $Y$ and $Z$, respectively, and 
$$T = {\sin{N}\over N}= \ \sum_{n\geq 0} {(-1)^n\over (2n+1)!}N^{2n}\ \in End(V)$$
Then
$$
z_b T y_a\not=0
$$
iff both $a$ and $b$ are weights, and  either $a$ is a highest weight  
or  $a<|b|$.
\endproclaim
\demo{Proof} Can take  $\rho$ irreducible of highest weight $d$ and realized
on the space of homogeneous polynomials of degree $d$ in two real variables 
$u$, $v$. The elements $u^{p}v^{d-p}$ form a basis adapted to $sl_2$, 
the elements $w^{p}\bar w^{d-p}$   ($w=u+iv$)  form a basis adapted 
to $su_2$ and  $N$ acts by $v{\partial \over \partial u}$. Since
 $$\eqalign{e^{iv{\partial \over \partial u}} u^p 
 &= \sum {1\over m!}({\partial \over \partial u})^m u^p(iv)^m  = \sum {p(p-1)\cdots (p-m+1)\over m!}u^{p-m}(iv)^m \cr
&= (u+iv)^p = w^p\cr}$$
one has
 $$
	\eqalign{e^{iv{\partial \over \partial u}} u^pv^{d-p} 
	&= w^pv^{d-p} 
	=  {1\over (2i)^{d-p}} w^p(w-\bar w)^{d-p},\cr 
e^{-iv{\partial \over \partial u}} u^pv^{d-p} 
	&=  {1\over (2i)^{d-p}} \bar w^p(w-\bar w)^{d-p}\cr}.
$$ Therefore
\goodbreak
$$
\eqalign{ {T} \ u^pv^{d-p}
&={\sin{N}\over {N}} u^pv^{d-p} \cr
	  &= (2i)^{-1} (e^{iv{\partial \over \partial u}}- 
   	       e^{-iv{\partial \over \partial u}})
	       (v{\partial \over \partial u})^{-1}\, u^pv^{d-p} \cr 
	  &= (2i)^{-1} (e^{iv{\partial \over \partial u}}- 
	       e^{-iv{\partial \over \partial u}})\,
	       {1\over p+1}u^{p+1}v^{d-p-1}\cr
	  &= {1\over (2i)^{d-p+1}(p+1)}
	       (w^{p+1}-\bar w^{p+1})(w-\bar w)^{d-p-1}\cr}
$$
For $p=d$ one gets
$$
 T\ u^d 
	  = {1\over 2i (d+1)}(w^{d+1}-\bar w^{d+1})(w-\bar w)^{-1} 	  = {1\over 2i (d+1)} \sum_{p=0}^d w^p\bar w^{d-p}
$$
which has non-zero projections on all $E^{Z}_j$'s. On the other hand, if 
$p<d$, then $d-p-1\geq 0$ and
$$
\eqalign{ (2i)^{d-p+1}(p+1){T} u^pv^{d-p} 
	  &= w^{p+1}(w-\bar w)^{d-p-1}-\bar w^{p+1}(w-\bar w)^{d-p-1}\cr
	  &= \sum_{j=0}^{d-p-1}{d-p-1\choose j}
	       \big(w^{d-j}\bar w^j- w^{j}\bar w^{d-j}\big)\cr}
$$
which has non-zero components exactly in
$ (E^{Z}_d \oplus \cdots \oplus E^{Z}_{2p-d+2})\oplus 
	(E^{Z}_{-(2p-d+2)}\oplus \cdots \oplus E^{Z}_{-d}).$
This finishes the proof of the Lemma.
\enddemo
%%%%%%%%%%%%%%%%%%%%%%%%%%%%% End section-4 %%%%%%%%%%%%%%%%%%%%%%%%%%%%%%%%
%%%%%%%%%%%%%%%%%%%%%%%%%%%%%% section-5 %%%%%%%%%%%%%%%%%%%%%%%%%%%%%%%%%%%

$$
	\text{\bf 5.\quad Split Variations}
$$

\par Here we will prove the following 
version of Theorem III:

\proclaim{Theorem 5.1} Let $\V\to\Delta^*$ be a variation of 
graded-polarized mixed Hodge structure with unipotent monodromy, and 
$\Y$ denote the $C^{\infty}$ grading the weight filtration $\W$ obtained by 
applying  $(1.2)$ to each fiber of $\V$.  Then, the variation $\V$ 
is admissible and split  if and 
only if
\roster
\item"(a)" $\Y$ extends continuously to a grading of $\W$ in the canonical 
prolongation of $\V$, with limiting value $\lim_{s\to 0}\, \Y(s)$ defined 
over $\R$.
\item"(b)" Each of the induced 
variations $\F Gr^{\W}_k$ is split.
\endroster
\endproclaim

As noted in the introduction, in order to make the
various reality conditions well defined, one must
select a non-zero reference element $\theta\in T_0(\Delta)^*$.

\par First, the
limiting Hodge filtration $F_{\infty}$ of such a variation always exists. If $\V$ is admissible, this is true by definition.  On the other hand, if (a) and (b) hold, then $F_{\infty}$ is obtained by extending the graded limiting Hodge filtrations to $V$ via $Y_{\infty}$, since $\Y$ preserves $\F$.
More explicitly, if $F:U\to\M$ represents the 
period map of $\V$ then
$$
	F_{\infty} = \lim_{s\to 0}\,\psi(s)			\tag{5.2}
$$
where $\psi(s):\Delta^*\to\check\M$ is the map associated to the periodic
function 
$$\psi(z) = e^{-zN}.F(z).$$ 
\goodbreak
Likewise, the hypothesis assert the
existence of the grading
$$
	Y_{\infty} = \lim_{\text{Im}(z)\to\infty}\, Y(z)
		   = \lim_{s\to 0}\, Y(s)			\tag{5.3}
$$
where $Y(s)$ is the grading of $W$ defined by the periodic function
$$Y(z) = e^{-zN}.Y_{(F(z),W)}.$$
Since  the graded variation is a sum of variations of pure Hodge structures,
$\lim_{s\to 0}\, \psi(s)Gr^W$ exists.  Since $Y(s)$ preserves the filtration $\psi(s)$, it follows that 
whenever $Y_{\infty}$ exists, so does $F_{\infty}$.
We may therefore write
$$
	F(z) = e^{zN}e^{\Gamma(s)}.F_{\infty},\qquad \Gamma(0) = 0
								\tag{5.4}
$$
relative to a suitable holomorphic function $\Gamma(s)$ (cf. $(2.3)$).  

To complete the proof, we will show that under the
hypothesis of the Theorem,
$$
	\lim_{\text{Im}(z)\to\infty}\, e^{-zN}.Y_{(F(z),W)}
	= \lim_{\text{Im}(z)\to\infty}\, e^{-zN}.Y_{(e^{zN}.F_{\infty},W)}
								\tag{5.5}
$$
Explanation: if  $(5.5)$ holds and $\V$ is admissible and split, $Y_s$ extends as described in
the introduction.  Conversely, if both $(5.5)$ and the hypothesis of the Theorem  hold, we can invoke $(4.1)$ to prove
that $\V$ is admissible, with limiting mixed Hodge structure 
$(F_{\infty},\rel W)$ which splits over $\R$.

\par Let $\Cal C$ to be the class of real-analytic 
functions $q(z)$ which take values in $GL(V)^W$ and satisfy an estimate 
of the form
$$
	q(z) = 1 + o(e^{-ky}),\qquad k>0
$$
for $Im(z)>>0$, and denote membership in $\Cal C$ by the shorthand $q\approx 1$.
The correspondence
$$
	(F,W) \leftrightarrow (Y_{(F,W)},F\,Gr^W)	 % was (3.3)
$$
is (real) birational and smooth along $\M$, since the $I^{p,q}$ are obtained from $F$ and $W$ by
taking finite intersections and sums.  Consequently, 
an $\M$-valued function $\phi(z)$ is of the form
$$
	\phi(z) = q(z).F_o
$$
relative to a fixed base point $F_o\in\M$ and some function $q(z)$ of 
class $\Cal C$ iff 
$$
	(Y_{(\phi(z),W)}, \phi(z) Gr^W)\rightarrow (Y_{(F_o,W)},F_o Gr^W)
$$ 
exponentially fast.  By smoothness, this conclusion remains valid for variable $F_o$, provided the limiting value
of $F_o$ belongs to $\M$.
Recall the notation $Y_F=Y_{(F,W)}$ for $F\in\M$.
\vskip 10pt

% Theorem (5.6) was (3.4)
\proclaim{Theorem 5.6} Let $F(z)$ be the period map of a split admissible variation
of graded-polarized mixed Hodge structure with unipotent monodromy .
Then, the limiting grading $(5.2)$ may be computed in terms of the 
corresponding nilpotent orbit, i.e.
$$
	\lim_{Im(z)\rightarrow\infty}\, e^{-zN}.Y_{F(z)} 
        = \lim_{Im(z)\rightarrow\infty}\, e^{-zN}.Y_{e^{zN}.F_{\infty}}
$$
\endproclaim
\demo{Proof} By $(5.3)$, we may write
our period map $F(z)$ in the form
$$
	F(z) = e^{zN}q(z).F_{\infty}
$$
relative to a function $q(z)\approx 1$.

\par Next, we note that splitting $\rel Y$ of $(F_{\infty},\rel W)$ preserves 
$F_{\infty}$ and satisfies $[\rel Y,N]=-2N$.  Consequently,
$$ 
	e^{iyN}.F_{\infty} =  y^{-\half\rel Y}.e^{iN}F_{\infty}
$$
Therefore, upon writing $q_1(z)=Ad(e^{iyN})q(z)$ and 
$q_2(z)=y^{\half \rel Y}.q_1(z)$, we obtain
$$
\align
	F(z) &= e^{xN}e^{iyN}q(z).F_{\infty}
	      = e^{xN}q_1(z)e^{iyN}.F_{\infty}				\\
	     &= e^{xN}q_1(z)y^{-\half\rel Y}e^{iN}.F_{\infty}
	      = e^{xN}y^{-\half\rel Y}q_2(z)e^{iN}.F_{\infty}		
\endalign
$$
Since $e^{iyN}$ and $y^{-\half\rel Y}$ preserve $W$ and are 
polynomials in positive and negative powers of $y^{\half}$, 
$q_1$ and $q_2\approx 1$.  Since $e^{iN}.F_{\infty}\in\M$ by virtue of Schmid's $SL_2$ Orbit Theorem:
$$
	Y_{F(z)}= Y_{e^{xN}y^{-\half\rel Y}q_2(z)e^{iN}.F_{\infty}}
	= e^{xN}y^{-\half\rel Y}.Y_{q_2(z)e^{iN}.F_{\infty}}
						\tag{5.7} % was (3.5)
$$

\par The map
$$
	\M\to Y(W),\qquad F\in\M \mapsto Y_{(F,W)}
$$
 into the space $Y(W)$ consisting of all
gradings $Y$ of $W$ is real-algebraic. Since $q_2\approx 1$, 
$e^{iN}.F_{\infty}\in \M$ and $\Gc$ acts transitively on $Y(W)$, we conclude:
$$
	Y_{q_2(z)e^{iN}.F_{\infty}} = q_3(z).Y_{e^{iN}.F_{\infty}}
$$
for some $q_3\approx 1$. Consequently, $(5.7)$ implies that
$$
\align
	  e^{-zN}.Y_{F(z)}
	  &= e^{-iyN}y^{-\half\rel Y}q_3(z).Y_{e^{iN}F_{\infty}}	\\
	  &= q_4(z)e^{-iyN}y^{-\half\rel Y}.Y_{e^{iN}.F_{\infty}}
	   = q_4(z)e^{-iyN}.Y_{e^{iyN}F_{\infty}}
\endalign
$$
for some $q_4\approx 1$, and hence:
$$
	\lim_{Im(z)\rightarrow\infty}\, e^{-zN}.Y_{F(z)} 
        = \lim_{Im(z)\rightarrow\infty}\, e^{-zN}.Y_{e^{zN}.F_{\infty}}
$$
\enddemo

\proclaim{Corollary} Under the hypothesis of Theorem $(5.6)$, the period
map $F(z)$ satisfies  $(a)$ and $(b)$.  Moreover, in this setting
the limiting grading $Y_{\infty}$ coincides with the grading of $W$ obtained by applying the first theorem of \cite{6} (cf\. Appendix) to the pair
$$
	(N,\rel Y),\qquad \rel Y = Y_{(F_{\infty},\rel W)}
$$
\endproclaim

\par To prove the converse, let $\rho$ be the representation of $sl_2(\R)$
attached to a variation $\V$ which satisfies  (a) and (b) of
Theorem (5.1), i.e. pull back the representations $\rho_k$ 
which define the associated nilpotent orbits $e^{zN}.FGr^W_k$ -which are
$SL_2$ Orbits, on account of  (b), via the grading 
$Y_{\infty}$, which is defined over $\R$.  Then, as in \S 4:
\roster
\item $N_0 = \rho(h)$, where
$$
	N = N_0 + N_{-1} + \cdots
$$
denotes the decomposition of $N$ with respect to the eigenvalues of 
$ad\,Y_{\infty}$.  

\item $Y_o = \rho(h)$ preserves $F_{\infty}$.
\endroster
Likewise, conditions (a) and (b) imply that $Y_{\infty}$
preserves $F_{\infty}$, since the grading
$Y(s)$ from  $(5.2)$ preserves the filtration $\psi(s)$.

% Theorem (5.8) was (3.10)
\proclaim{Theorem 5.8}  Let $F(z)$ be the period map of a variation of 
graded-polarized mixed Hodge structure with unipotent monodromy for 
which conditions (a) and (b) of Theorem (5.1) hold.  Then, 
$$
	Y_{\infty} = \lim_{Im(z)\to\infty}\, e^{-zN}.Y_{e^{zN}.F_{\infty}}
$$
\endproclaim
\demo{Proof} By $(5.3)$, 
$$
	F(z) = e^{zN}g(z).F_{\infty}
$$
for some function $g(z)\approx 1$.  Since $\Cal C$ is
closed under conjugation by $e^{zN}$, 
$$
	F(z)=q(z)e^{zN}.F_{\infty}
$$
for some function $q\approx 1$.
Next, we prove that for a suitable integer $m\in\Z$, the endomorphism
$$
	T = Y_o + m Y_{\infty},\qquad Y_o = \rho(H)
$$
satisfies the estimate
$$
	yAd(y^T)N = N_0 + o(y^{-1})
$$
Indeed, decompose $N=N_0+N_{-1}+\cdots$ according to the eigenvalues of
$Y_{\infty}$.  Then,
$$
	yAd(y^T)N =  Ad(y^{\half Y_o}) y Ad(y^{mY_{\infty}})N  
	= N_0 + \sum_{j\geq 1} y^{1-mj}Ad(y^{\half Y_o})N_{-j}.
$$
Since $j\geq 1$ and only finitely many 
eigenvalues of $Y_o$ appear non-trivially in each $N_j$, we may force the leftmost sum to involve only negative powers of $y$, by making $m$ sufficiently large .

\par  By hypothesis, $T$ is defined over $\R$. The remarks immediately
preceding $(5.8)$, $T$ also preserves 
$F_{\infty}$. Hence
$$
	e^{iyN}.F_{\infty} 
	= e^{iyN}y^{-T}.F_{\infty} 
	= y^{-T}e^{iyAd(y^T)N}.F_{\infty} 
	= y^{-T}P(y)e^{iN_0}.F_{\infty}
$$
relative to a function $P(y)$ which is polynomial in $y^{-1}$.

\par To finish the proof, we note that by virtue of  $(b)$ and
the $SL_2$ Orbit Theorem, $e^{iN_0}.F_{\infty}\in \M$.
Since $\Cal C$ is closed under conjugation by 
$y^{kT}$, $k\in\R$:
$$
\aligned
	 F(iy)  &=  q(iy)e^{iyN}.F_{\infty}  
	         = q(iy) y^{-T} e^{iyAd(y^T)N}.F_{\infty}	\\ 
		&= y^{-T} q_1(y) e^{iyAd(y^T)N}.F_{\infty}
\endaligned
$$
and hence
$$
	Y_{F(iy)} = y^{-T}.Y_{q_1(iy) e^{iyAd(y^T)N}.F_{\infty}}
		  = y^{-T} q_2(y).Y_{e^{iyAd(y^T)N}.F_{\infty}}.
$$
Accordingly,
$$
\eqalign{
	e^{-iyN}.Y_{F(iy)} 
	&= e^{-iyN}y^{-T}q_2(y).Y_{e^{iyAd(y^T)N}.F_{\infty}}  
	 = q_3(y)e^{-iyN}y^{-T}.Y_{e^{iyAd(y^T)N}.F_{\infty}} \cr 
	&= q_3(y)e^{-iyN}.Y_{y^{-T}e^{iyAd(y^T)N}.F_{\infty}}  
	 = q_3(y)e^{-iyN}.Y_{e^{iyN}y^{-T}.F_{\infty}} \cr 
	&= q_3(y)e^{-iyN}.Y_{e^{iyN}.F_{\infty}}} 
$$
for functions $q_2$ and $q_3$ $\approx 1$, wherefrom
$$
	\lim_{y\to\infty}\, 
	e^{-iyN}.Y_{e^{iyN}.F_{\infty}} = Y_{\infty}
$$
\enddemo

%%%%%%%%%%%%%%%%%%%%%%%%%%%% End section-5 %%%%%%%%%%%%%%%%%%%%%%%%%%%%%

%%%%%%%%%%%%%%%%%%%%%%%%%%%%%% appendix %%%%%%%%%%%%%%%%%%%%%%%%%%%%%%%%

$$
	\text{\bf Appendix.\quad Two theorems of Deligne}
$$

\par Let 
$ 0\subseteq\cdots\subseteq W_k \subseteq W_{k+1}\subseteq\cdots\subseteq V$
be an increasing filtration of a finite dimensional vector space $V$ defined
over an algebraically closed field $k$ of characteristic zero.
\vskip 10pt

\proclaim{Theorem 1} Let $N$ be a nilpotent endomorphism of $V$ for which 
the corresponding relative weight filtration $\rel W = \rel W(N,W)$ exists,
and $\rel Y$ be a grading of $\rel W$ which preserves $W$ and satisfies the
additional condition $[\rel Y, N] = -2N$. Then, there exists a unique grading 
$Y$ of $W$ such that
\roster
\item $[\rel Y,Y] = 0$.
\item $N = N_0 + N_{-2} + \cdots$ when decomposed relative to the
      eigenvalues of $ad\,Y$.
\item $(ad N_0)^{k-1}\,N_{-k} = 0$ for all $k>0$.
\endroster
\endproclaim
\demo{Sketch of Proof} The desired grading $Y$ may be constructed as
follows: Let $Y_o$ be a grading of $W$ for which $[\rel Y,Y] = 0$, and 
note that the group 
$$
        G_o = \{\, g\in GL(V) \mid [g,\rel Y] = 0,\quad 
                                   (g-1):W_k\to W_{k-1}\hph{a}\forall k \,\}
$$
acts transitively on the set of all such gradings.  Assume by induction
that the initial grading $Y_o$ satisfies  $(1)$--$(3)$ modulo the 
ideal
$$
	Lie_{-r} = \{\, \a\in End(V) \mid \a:W_k\to W_{k-r}\hph{a}\forall k\,\}
$$
Then one may construct an element 
$g\in G_o$ such that, relative to the grading $Y' = Ad(g)Y_o$,  
$(1)$--$(3)$ hold modulo the ideal $Lie_{-(r+1)}\subseteq Lie_{-r}$.
Since $Lie_{-r} = 0$ for some finite index $r$, the desired grading $Y$
will be obtained after finitely many steps. The details 
 are discussed in \cite{12}.
\enddemo

\flushpar{\it Remark.} The preceding result may be reformulated as the 
statement that given the existence of $\rel W = \rel W(N,W)$ and a grading 
$\rel Y$ of the type described above, there exists a unique grading $Y$ of $W$
which commutes with $\rel Y$ and has the property that the 
associated $sl_2$-triple $(N_0,\rel Y-Y,N_0^+)$ satisfies the commutativity 
condition $[N - N_0,N_0^+] = 0$.  

\proclaim{Lemma} The construction of Theorem 1 is both functorial
and compatible with the operations of direct sum, tensor product and dual.  
\endproclaim

\par To state Deligne's second theorem, we begin with a triple $(F,W,N)$ of 
the type arising from the degeneration of an admissible variation of graded-polarizable 
mixed Hodge structure, i.e. 
\roster
\itemb $V = V_{\R}\otimes\C$ relative to some underlying real 
       form $V_{\R}$ to which both $W$ and $N$ descend.
\itemb The relative weight filtration $\rel W = \rel W(N,W)$ 
       exists.
\itemb The pair $(F,\rel W)$ is a mixed Hodge structure relative to which $N$ is a $(-1,-1)$-morphism and $W$ is a filtration by sub-mixed Hodge
       structures.
\endroster
In particular, given any such triple $(F,W,N)$, we may construct an
associated grading 
$$
	Y = Y(F,W,N)					
$$
of $W$ via application of Theorem 1 to the grading $\rel Y$ of $\rel W$ 
determined by the $I^{p,q}$'s of the mixed Hodge structure $(F,\rel W)$.
\vskip 10pt

\flushpar{\it Remark.} Inspection of the proof of Theorem 1 shows that 
whenever the mixed Hodge structure $(F,\rel W)$ associated to the triple $(F,W,N)$ is split, the corresponding grading $Y=Y(F,W,N)$ lies in
$End(V_{\R})$.
\vskip 10pt

\proclaim{Lemma} The grading $Y=Y(F,W,N)$ constructed above preserves $F$.
\endproclaim
\demo{Sketch of Proof} (cf.  \cite{12} for details).  One uses the preceding lemma to reduce the problem to the case where 
$(F,\rel W)$ is a split mixed Hodge structure which is either of the form
\roster
\itemb $V = \oplus_p\, I^{p,p}$ with $N(I^{p,p}) \subseteq I^{p-1,p-1}$ or, 
\itemb $V = I^{p,0}\oplus I^{0,p}$ and $N=0$.
\endroster
\enddemo

%
% Double check the statement of Theorem 2, make sure we don't need any
% extra conditions [for example: graded-polarizations may be needed to
% ensure that $(e^{iN}.F,\rel W)$ is a MHS.]
%

\proclaim{Theorem 2} Let $(F,W,N)$ be a triple of the type described
above for which the  mixed Hodge structure $(F,\rel W)$ is split.  Then, the pair $(e^{iyN}.F,\rel W)$ defines a mixed Hodge
structure for all $y>0$.  Moreover, the grading $Y_{(e^{iyN}.F,W)}$ of $W$ 
determined by the $I^{p,q}$'s of $(e^{iyN}.F, W)$ is given by the formula
$$
	Y_{(e^{iyN}.F,W)} = e^{iyN}.Y(F,W,N)			
$$
\endproclaim
\demo{Proof} As a consequence of \cite{3}, the pair $(e^{iyN}.F, W)$ 
is a mixed Hodge structure for all $y>0$.  To establish the stated formula, we define $Y$ to
be $Y(F,W,N)$, and note that
$$
	Y = Y_{(e^{iyN_0}.F,W)}  				
$$
for all $y>0$ since, in the present context, $Y$ is a grading of $W$ which 
is both defined over $\R$ and preserves $e^{iyN_0}.F$.  Accordingly, in 
order to establish Theorem 2, it will suffice to show that
$$
	e^{\phi(y)} := e^{iyN}e^{-iyN_0}\in\exp(\lam_{(e^{iyN_0}.F,W)})     
								\tag{\ast}
$$
Indeed, if the preceding equation is true, then 
$$
\aligned
	Y_{(e^{iyN}.F,W)} 
	 &= Y_{(e^{iyN}e^{-iyN_0}e^{iyN_0}.F,W)}			\\
	 &= e^{iyN}e^{-iyN_0}.Y_{(e^{iyN_0}.F,W)}			
	 = e^{iyN}e^{-iyN_0}.Y
	 = e^{iyN}.Y	
\endaligned
$$
\par To verify  $(\ast)$, note that $\phi(y)$ is a Lie polynomial in the 
monomials $(ad\,N_0)^a\,N_{-k}$.  In particular $\phi(y)$ contains no $N_0$ 
term.  Moreover, by conditions $(2)$ and $(3)$ of Theorem 1, 
$(ad\,N_0)^a\,N_{-k}$ is zero unless $k>1$ and $0\leq a\leq k-2$.  It then follows from $(4.10)$  that 
$$
	(ad\,N_0)^a\,N_{-k} 
	\in \bigoplus_{0< r <k}\, gl(V)^{-r,r-k}_{(e^{iyN_0},W)}
	\subseteq\lam_{(e^{iyN_0}.F,W)}
$$
Consequently, $\phi(y)$ takes values in $\lam_{(e^{iyN_0}.F,W)}$ since the
latter is a Lie subalgebra of $gl(V)$.
\enddemo

\

It is clear that Theorem (IV) in the introduction follows as a corollary of these proofs.
%%%%%%%%%%%%%%%%%%%%%%%%%%%%%%% End appendix %%%%%%%%%%%%%%%%%%%%%%%%%%%%

%%%%%%%%%%%%%%%%%%%%%%%%%%%%%%%% References %%%%%%%%%%%%%%%%%%%%%%%%%%%%%%


\Refs

% \widestnumber


%\key BZ
\ref \no 1 \by Brylinski J., and Zucker S.
\paper An Overview of Recent Advances in Hodge Theory.
\inbook Complex Manifolds, (Gamkrelidze R. ed.), Springer 
\yr 1997
\endref

%\key CKS
\ref \no 2 \by Cattani E., Deligne P.,  and Kaplan A. 
\paper On the locus of Hodge cycles
\jour Journal Amer. Math. Soc. \vol 8 \pages 483--506 \yr 1995 
\endref

%\key CG

\ref \no 3 \by Cattani E., Kaplan A.,  and Schmid W. 
\paper Degeneration of Hodge structures
\jour Ann. of Math. \vol 123 \pages 457--535 \yr 1986
\endref

\ref \no 4 \by Cattani E., Kaplan A. and Schmid W. 
\paper $L_2$ and Intersection cohomologies with coefficients in a variation of Hodge structure.
\jour Invent. Math. \vol 87 \pages 217--252 \yr 1986
\endref

\ref \no 5 \by Cornalba M. and Griffiths P.
\paper Some transcendental aspects of algebraic geometry
\jour Proceedings of Symposia in Pure Mathematics
\vol 29 \pages 3--110 \yr 1975
\endref
%\key D1
\ref \no 6 \by Deligne P.
\paper Private communication \yr 1995
\endref


%\key D2
\ref \no 7 \by Deligne P.
\paper {\rm Appendix to} Variation of Mixed Hodge Structure I
\jour Invent. Math. \vol 80, \pages 489 -- 542 \yr 1985
\endref

%\key D3
% Note: Was D4
\ref \no 8 \by Deligne P.
\paper Equations diff\'erentielles \`a singuliers r\'eguliers
\inbook  Lecture Notes in Mathematics  \vol 163  \yr 1970.
\endref

%\key HZ
\ref \no 9 \by Hain R., and Zucker S. 
\paper Unipotent variations of mixed Hodge structure
\jour Invent. Math. \vol 88, \pages 83 -- 124 \yr 1987
\endref

% \key K
\ref \no 10 \by Kaplan A.
\paper Notes on the moduli spaces of Hodge structures
\jour Unpublished manuscript \yr 1995 
\endref

\ref \no 11 \by Kato, K., Usui, S. 
\paper Logarithmic Hodge structures and classifying spaces
\jour CRM Proc. Lecture Notes, Amer. Math. Soc.\vol 24\pages 115--130 \yr 2000 
\endref

%\key P1
\ref \no 12 \by Pearlstein G.
\paper Degenerations of Mixed Hodge Structure, [math.AG/0002030]  
\inbook To appear, Duke Math. Journal 
\endref

\ref \no 13 \by Pearlstein G.
\paper $SL_2$ Orbits and Degenerations of Mixed Hodge Structure
\inbook Preprint
\endref

%\key P2
\ref \no 14 \by Pearlstein G.
\paper Variations of Mixed Hodge Structure, Higgs Fields and Quantum 
Cohomology, [math.AG/9808106]
\jour Manuscripta Math. \vol 102 \yr 2000 \pages 269--310
\endref

%\key S
\ref \no 15 \by Schmid W.
\paper Variation of Hodge Structure: The Singularities of the Period Mapping.
\jour Invent. Math. \vol 22 \pages 211 -- 319  \yr 1973
\endref

%\key SZ
\ref \no 16 \by Steenbrink J., and Zucker S.
\paper Variation of Mixed Hodge Structure I.
\jour Invent. Math. \vol 80 \pages 489 -- 542 \yr 1985
\endref
\endRefs
\enddocument